\documentclass{amsart}

\usepackage{etoolbox} 

 \usepackage{amssymb}
\usepackage[pdfauthor="Richard J. Mathar",colorlinks=true,bookmarksopen,pdfpagemode=UseOutlines]{hyperref}

\newtheorem{defn}{Definition}
\newtheorem{rem}{Remark}

\usepackage{orcidlink}
\begin{document}

\title[Chebyshev Expansion of Zernike Polynomials]{Chebyshev Expansion of the Zernike Polynomials}

\author{Richard J. Mathar \orcidlink{0000-0001-6017-6540}}
\urladdr{https://mathar.www3.mpia.de}
\email{mathar@mpia.de}
\address{Max-Planck Institute for Astronomy, K\"onigstuhl 17, 69117 Heidelberg, Germany}


\date{\today}

\begin{abstract}
The even  and odd Zernike Polynomials $R_n^m(x)$
can be expanded into sums of even and odd Chebyshev Polynomials $T_i(x)$.
This manuscript provides closed forms for the rational expansion coefficients 
for a set of small $0\le n-m\le 6$ and a holonomic five-term recurrence 
for these coefficients
for all larger $n$.
\end{abstract}

\maketitle

\section{Overview}
Zernike Circle Polynomials $R_n^m(x)$ are polynomials of degree $n$
in the radial coordinate $x$, $0\le x \le 1$
and defined for even non-negative $n-m=0,2,4,\ldots$, $0\le m\le n$:
\begin{defn} (2D Zernike Radial Polynomials)
\begin{multline}
R_n^m(x)
=
(-1)^{(n-m)/2}
x^m\sum_{s=0}^{(n-m)/2}
\binom{(n+m)/2+s}{(n-m)/2-s}\binom{m+2s}{s}(-x^2)^s
\\
=
(-1)^{(n-m)/2}
\binom{(n+m)/2}{(n-m)/2}
x^m
{}_2F_1\left(\begin{array}{c} -(n-m)/2, 1+(n+m)/2 \\ m+1\end{array}\mid x^2\right)
\label{eq.Rofx}
\end{multline}
where $_pF_q$ are generalized hypergeometric series.
\end{defn}
We shall call $n$ the radial quantum number and $m$ the azimuthal
quantum number.
The inverse relation expands monomials $x ^j$ for fixed $m$ into $R_n^m$
by defining rational coefficients $h$ \cite{MatharSAJ179},
where $j-m=0,2,4,\ldots$.

The Chebyshev Polynomials $T_i(x)$ are another set of polynomials
defined over $-1\le x\le 1$:
\begin{defn} (Chebyshev Polynomials of the First Kind) \cite[22.3]{AS}
\begin{equation}
T_0(x)=1;\quad
T_i = 
\frac{i}{2}\sum_{s=0}^{\lfloor i/2\rfloor}
(-)^s \frac{(i-s-1)!}{s!(i-2s)!}(2x)^{i-2s},\quad i=1,2,3,4\ldots
\label{eq.T}
\end{equation}
\end{defn}
The are orthogonal with an inverse square root weight:
\begin{equation}
\int_{-1}^1 T_j(x) T_k(x)\frac{dx}{\sqrt{1-x^2}}= \delta_{j,k} \frac{\pi}{2} \epsilon_j,\quad j,k\ge 0.
\label{eq.Tortho}
\end{equation}
\begin{defn} (Neumann symbol)
\begin{equation}
\epsilon_n =\left\{\begin{array}{ll} 2, & n=0; \\ 1 & n>0.  \end{array}\right.
\end{equation}
\end{defn}
The reverse expansion of \eqref{eq.T} is
\cite[p.\ 412]{Cody}\cite[p.\ 52]{Fox}\cite{FraserJACM12}
\begin{equation}
x^n=2^{1-n}{\mathop{{\sum}'}_{j=0 \atop n-j\,\mathrm{even}}^n}{n\choose (n-j)/2} T_j(x) ,
\label{eq.xnofT}
\end{equation}
where the prime at the sum symbol means the term at summation index $j=0$
is to be halved if present.

The expansion of the $T_i$ in finite series of $R_n^m$ 
is known \cite{MatharVixra2403}; this manuscript derives
the expansion coefficients of the inverse problem. 
The main impetus to seek a re-expansion of the Zernike polynomials in other polynomials is 
that a summation of alternating monomials derived straight from \eqref{eq.Rofx} becomes
numerically unstable for large $n$.
Aside from other methods to compute $R_n^m(x)$ by
recurrences 
\cite{KintnerJModOpt23,ChongPR36},
this is an alternative of approaching the problem of the cancellation of digits 
(which happens for large $n$ while accumulating the alternating terms of the monomial basis) 
\cite{HoudayerHDMIalg15}.

\begin{defn}(Coupling coefficients)\label{def.c}
\begin{equation}
R_n^m(x) \equiv \mathop{{\sum}'}_i
c_{n,m,i} T_i(x)
\label{eq.cdef}
\end{equation}
for even $n-m \ge 0$.
\end{defn}

\begin{rem}
Both types of polynomials obey second order differential equations,
so another approach---not followed here---is to start from the differential
equation of the Zernike polynomials \cite{MatharArxiv0705a} and to use the generic algorithms
of seeking solutions in the Chebyshev basis \cite{ClenshawMPCPS53}.
\end{rem}

A basic result of comparing the parities of the polynomials in \eqref{eq.cdef} 
are the selection rules
\begin{equation}
c_{n,m,i}=0, \mathrm{ if} n-m \nmid 0 \pmod 2 \vee n-i \nmid 0 \pmod 2.
\label{eq.selec}
\end{equation}
With a computer algebra system at hand, one can compute the $c_{n,m,i}$ recursively
by a basic form of long polynomial division and deflation:
Starting with $i=n$, $c_{n,m,i}$ is given by the ratio of the leading coefficient of
$R_n^m$ divided by the leading coefficient of $T_i$. A residual $R_n^m$ is obtained
by subtracting $c_{n,m,i}T_i(x)$ from $R_n^m(x)$, and the next lower coefficient $c_{n,m,i-2}$
is the ratio of the leading coefficient of the residual $R_n^m(x)$ divided by the leading coefficient
of $T_{i-2}(x)$---with the exception at the constant term where the prime at the sum
symbol in \eqref{eq.cdef} 
means $c_{m,n,0}$ is twice the ratio. This loop is executed until the residual $R_n^m$ is zero.
This produces Table \ref{tab.c}.

\begin{table}
\begin{tabular}{rr|rrrrrr}
$n$ & $m$ & $c_{n,m,i}$ \\
\hline
0& 0& 2 \\
1& 1& 1 \\
2& 0& 0& 1 \\
2& 2& 1& 1/2 \\
3& 1& 1/4& 3/4 \\
3& 3& 3/4& 1/4 \\
4& 0& 1/2& 0& 3/4 \\
4& 2& 0& 1/2& 1/2 \\
4& 4& 3/4& 1/2& 1/8 \\
5& 1& 1/4& 1/8& 5/8 \\
5& 3& 1/8& 9/16& 5/16 \\
5& 5& 5/8& 5/16& 1/16 \\
6& 0& 0& 3/8& 0& 5/8 \\
6& 2& 3/8& 1/32& 5/16& 15/32 \\
6& 4& 0& 5/16& 1/2& 3/16 \\
6& 6& 5/8& 15/32& 3/16& 1/32 \\
7& 1& 9/64& 15/64& 5/64& 35/64 \\
7& 3& 15/64& 1/64& 27/64& 21/64 \\
7& 5& 5/64& 27/64& 25/64& 7/64 \\
7& 7& 35/64& 21/64& 7/64& 1/64 \\
8& 0& 9/32& 0& 5/16& 0& 35/64 \\
8& 2& 0& 9/32& 1/16& 7/32& 7/16 \\
8& 4& 5/16& 1/16& 1/8& 7/16& 7/32 \\
8& 6& 0& 7/32& 7/16& 9/32& 1/16 \\
8& 8& 35/64& 7/16& 7/32& 1/16& 1/128
\end{tabular}
\caption{Reference values of $c_{n,m,i}$ for all $n\le 8$. The $c$-values 
are sorted left-to-right in each line by increasing $i=0,2,4,\ldots,n$ if $n$ is even, by
increasing $i=1,3,5,7,\ldots,n$ if $n$ is odd.
}
\label{tab.c}
\end{table}

Heuristically tested for $n\le 200$, all $c_{n,m,i}$ are positive---I have no formal proof of this---and
in conjunction with the sum rule \eqref{eq.srule} the values (except $c_{0,0,0}$) are $\le 1$.
Since $|T_i(x)|\le 1$, cancellation of digits is limited in these expansions to the
oscillations of the Chebyshev polynomials.

\section{Coupling via Monomials}\label{sec.M}
The strategy to compute the $c$ is to insert \eqref{eq.xnofT} into
the right hand side of \eqref{eq.Rofx}:
\begin{multline}
R_n^m(x)
=(-)^{(n-m)/2}
\sum_{s=0}^{(n-m)/2} (-)^s
\binom{(n+m)/2+s}{(n-m)/2-s}\binom{m+2s}{s}x^{m+2s}
\\
=(-)^{(n-m)/2}
\sum_{s=0}^{(n-m)/2} (-)^s
\binom{(n+m)/2+s}{(n-m)/2-s}\binom{m+2s}{s}
\\ \times
2^{1-m-2s}{\mathop{{\sum}'}_{i=0 \atop m+2s-i\,\mathrm{even}}^{m+2s}}{m+2s\choose (m+2s-i)/2} T_i(x) ,
\end{multline}
which in view of Defn.\ \ref{def.c} yields
\begin{multline}
c_{n,m,i}
=(-)^{(n-m)/2}
\sum_{s=\max(0,(i-m)/2)}^{(n-m)/2} (-)^s
\binom{(n+m)/2+s}{(n-m)/2-s}\binom{m+2s}{s}
\\ \times
2^{1-m-2s}
{m+2s\choose (m+2s-i)/2}
\label{eq.csum}
\end{multline}
for $0\le i \le n$ and even $m-i$.
\subsection{The case $i\le m$}
If $i\le m$, the lower $s$-limit in \eqref{eq.csum} is always zero,
\begin{multline}
c_{n,m,i}
=(-)^{(n-m)/2}
\sum_{s=0}^{(n-m)/2} (-)^s
\binom{(n+m)/2+s}{(n-m)/2-s}\binom{m+2s}{s}
\\ \times
2^{1-m-2s}
{m+2s\choose (m+2s-i)/2}
\end{multline}
and conversion of the binomial coefficients to $\Gamma$-function ratios
provides a representation as a terminating generalized hypergeometric series \cite{SlaterHyp,RoyAMM94}:
\begin{multline}
c_{n,m,i}=(-)^{(n-m)/2}
2^{1-m}
\frac{(\frac{n+m}{2})!}{(\frac{n-m}{2})!(\frac{m-i}{2})!(\frac{m+i}{2})!}
\\ \times
{}_4F_3\left(\begin{array}{c}
(m+1)/2, -(n-m)/2, 1+(n+m)/2, m/2+1\\
m+1, 1+(m-i)/2, 1+(m+i)/2
\end{array}\mid 1\right)
.
\label{eq.4F3}
\end{multline}
The parametric excess (sum of the three lower parameters minus sum of the four
upper parameters) is $1/2$.

A table at small nonnegative 
\begin{equation}
\kappa\equiv (n-m)/2,
\end{equation} 
$i\le m$, can  be organized 
as
\begin{equation}
c_{n,m,i} 
=
c_{m+2\kappa,m,i} 
= 
2^{1-m}
\binom{ (n+m)/2}{\kappa}
\binom{m}{(m-i)/2}
\frac{\omega_{\kappa}(m,i) }
{\prod_{d=1}^\kappa (2d+m+i)(2d+m-i)},\quad i\le m.
\label{eq.m2m}
\end{equation}
Empty products (upper limit smaller than the lower limit)
equal 1. 

With a notational shortcut
\begin{defn}
\begin{equation}
m_+\equiv (m+i)/2; \quad m_- \equiv (m-i)/2, 
\end{equation}
\end{defn}
this is also
\begin{equation}
c_{n,m,i} 
= 
2^{1-m-2\kappa}
\binom{ (n+m)/2}{\kappa}
\binom{m}{m_-}
\frac{\omega_{\kappa}(m,i) }
{\prod_{d=1}^\kappa (d+m_+)(d+m_-)},\quad i\le m.
\end{equation}
For $n-m\le 6$ the numerators in this formula are in Table \ref{tab.omeg}.
\begin{table}
\begin{tabular}{r|c}
$\kappa$ & $\omega_\kappa(m,i)$\\
\hline
0 & 1 \\
1 & $i^2$ \\
2 & $(-m+i^2-4)^2$ \\
3 & $i^2(-3m-16 +i^2)^2$ \\
4 & $(3m^2+42m+144-6mi^2-40i^2+i^4)^2$ \\
5 & $i^2(15m^2+250m+1024-10mi^2-80i^2+i^4)^2$ \\
6 & $(-15m^3-450m^2-4440m-14400+45m^2i^2+870mi^2+4144i^2-15mi^4-140i^4+i^6)^2$
\end{tabular}
\caption{Table of numerators $\omega$ in \eqref{eq.m2m} for small indices $\kappa$.}
\label{tab.omeg}
\end{table}

Another simpler case is when $i=0$---which implies $m$ and therefore $n$ are 
even---, because then the upper parameter $1+m/2$ and lower $1+(m+i)/2$ parameter
of the hypergeometric series \eqref{eq.4F3} cancel
and the result is the square of a Gaussian Hypergeometric Function
\cite[Entry 13]{BerndtBLMS15}\cite{ClausenCrelle3}
which has a well-known value for argument 1 \cite[15.1.20]{AS}
\begin{multline}
c_{n,m,0}
=(-)^{(n-m)/2}
2^{1-m}
\frac{(\frac{n+m}{2})!}{(\frac{n-m}{2})!(\frac{m}{2})!^2}
\\ \times
{}_3F_2\left(\begin{array}{c}
(m+1)/2, -(n-m)/2, 1+(n+m)/2\\
m+1, 1+m/2
\end{array}\mid 1\right)
\\
=(-)^{(n-m)/2}
2^{1-m}
\frac{(\frac{n+m}{2})!}{(\frac{n-m}{2})!(\frac{m}{2})!^2}
{}_2F_1^2\left(\begin{array}{c}
-(n-m)/4, 1/2+(n+m)/4\\
m/2+1
\end{array}\mid 1\right)
\\
= 
\begin{cases}
2^{1-(n+m)/2}
 \frac{(\frac{n+m}{2})!}{(\frac{n-m}{2})!}
\left(
\frac{[(n-m)/2-1]!! }{
((n+m)/4)!
 }
\right)^2, & (n-m)/2\, \textrm{even}; \\
0, & (n-m)/2\, \textrm{odd}.
\end{cases}
\label{eq.cnm0}
\end{multline}

Another simplifying case is $i=1$---which implies $m$ and therefore $n$ are odd---, because
the upper parameter $(1+m)/2$ and lower parameter $1+(m-i)/2$ in \eqref{eq.4F3} cancel:
\begin{multline}
c_{n,m,1}
=(-)^{(n-m)/2}
2^{1-m}
\frac{(\frac{n+m}{2})!}{(\frac{n-m}{2})!(\frac{m-1}{2})!(\frac{m+1}{2})!}
\\ \times
{}_3F_2\left(\begin{array}{c}
-(n-m)/2, 1+(n+m)/2, m/2+1\\
m+1, 1+(m+1)/2
\end{array}\mid 1\right)
.
\label{eq.cnm1}
\end{multline}

This matches the 3-term recurrence \cite[1.4]{BaileyPGMA2}
\begin{multline}
(2\sigma-k)(\sigma-k)(\beta+\delta-k-1)(\alpha+\beta-k-1)L_{k+2}
+Q_kL_{k+1}
\\
+(k+1)(\sigma-k-1)(\gamma+\delta-k)(\alpha+\gamma-k)L_k=0
\end{multline}
where 
\begin{equation}
L_k\equiv {}_3F_2(-k,-\beta,k-2\sigma-1 ; -\beta-\delta,-\alpha-\beta ;1),
\end{equation}
\begin{equation}
Q_k=(2\sigma-2k-1)[(\beta\gamma-\alpha\delta)(\sigma+1)+\theta_3(k+1)(2\sigma-k)],
\end{equation}
\begin{equation}
2\sigma=\alpha+\beta+\gamma+\delta,
\end{equation}
\begin{equation}
2\theta_3=\alpha-\beta-\gamma+\delta.
\end{equation}
If we set

$$
\alpha =-1/2;\,
\sigma = -m/2-1 ;\,
\delta=-m/2;\,
\beta=-m/2-1;\,
k=(n-m)/2;
$$
the explicit write-up is
\begin{multline}
\frac{(4+m+n)^2(2+n)(5+n)}{16} 
\\
\times
{}_3F_2\left(\begin{array}{c} 
-(n+4-m)/2 , 1+(n+4+m)/2, 1+m/2 \\
1+m, 3/2+m/2
\end{array}\mid 1\right)
\\
\frac{(3+n)}{8}(6n+n^2-m^2+8)
{}_3F_2\left(\begin{array}{c} -(n+2-m)/2 , 1+(n+2+m)/2, 1+m/2
\\ 
1+m, 3/2+m/2
\end{array}\mid 1\right)
\\
-\frac{(n-m+2)^2(4+n)(1+n)}{16}
{}_3F_2\left(\begin{array}{c} 
-(n-m)/2 , 1+(n+m)/2, 1+m/2
\\ 
1+m, 3/2+m/2
\end{array}\mid 1\right)
=0.
\end{multline}

Switching back from the three $_3F_2$ 
to the $c$-parameters with \eqref{eq.4F3}
and shifting $n\to n-4$ yields
\begin{multline}
-\left(\frac{n-m}{2}-1\right)\left(\frac{n+m}{2}-1\right)\frac{n-3}{2} nc_{n-4,m,1}
\\
-\left( \frac{m+n}{2}\frac{n-m}{2} -\frac{n}{2}\right) (n-1) c_{n-2,m,1}
\\
+\frac{m+n}{2}\frac{n-m}{2}\left(\frac{n}{2}-1\right) (n+1)c_{n,m,1} =0
.
\label{eq.rec1}
\end{multline}
So a strategy for $i=1$ is to take \eqref{eq.mm} or \eqref{eq.m2m}
if $n\le m+2$, and to construct the values for $n>m+2$ recursively
from smaller $n$ with the 2-term reductions of \eqref{eq.rec1}.

Two more special cases:
\begin{itemize}
\item
If $m=0$,
the upper parameter $1+m/2$ cancels the lower parameter $1+m$
in \eqref{eq.4F3} and only a $_3F_2$ function needs to be calculated.
\item
If $m=1$, the upper parameter $(m+1)/2$ is 1 and an inhomogeneous
recurrence relation of order 2 can be derived \cite{LewanowiczMC45}.
This is not executed in practice, because $i\le m$ requires $i=1$,
the lower parameter $1+(m-i)/2$ is also 1, so it cancels the upper parameter;
so only a $_3F_2$ function needs to be calculated.
\end{itemize}

Fields, Luke and Wimp published
a 5-term recurrence 
\begin{equation}
P_k(m) + \sum_{s=1}^4 [A_s(k,m)+B_s(k,m)]P_{k-s}(m)=0
\label{eq.5rec}
\end{equation}
for the generalized hypergeometric polynomials
\eqref{eq.4F3} \cite[(2.44)]{FieldsJAT1}\cite{WimpMathComp22}
\begin{equation}
P_k(m)\equiv 
{}_4F_3\left(\begin{array}{c}
-k,k+m,(m+1)/2,1+m/2\\
m+1,1+(m-i)/2,1+(m+i)/2
\end{array}\mid 1\right)
.
\end{equation}
We compute the four $A$ and $B$ with 
argument $z=1$, upper parameters
$\lambda=1+m$,
$\alpha_1=(m+1)/2$, and
$\alpha_2=1+m/2$,
and lower parameters
$\beta_1=m+1$,
$\beta_2=1+(m-i)/2$,
$\beta_3=1+(m+i)/2$.
Then substituting $k=(n-m)/2$ our particular case of the recurrence reads
\begin{multline}
{}_4F_3\left(\begin{array}{c}
-(n-m)/2,1+(n+m)/2,(m+1)/2,1+m/2\\
m+1,1+(m-i)/2,1+(m+i)/2
\end{array}\mid 1\right)
\\
+
\frac{4(n-1)(4n-4n^2-m^2n-ni^2+n^3+m^2i^2)}{(n+i)(n-i)(n-4)(n+m)^2}
\\
\times
{}_4F_3\left(\begin{array}{c}
-(n-m)/2+1,1+(n+m)/2-1,(m+1)/2,1+m/2\\
m+1,1+(m-i)/2,1+(m+i)/2
\end{array}\mid 1\right)
\\
-2
(n-m-2)(n-3)n
\\ \times
\frac{32-72n-2m^2-3m^2i^2-2i^2+48n^2+6m^2n-n^2i^2+6ni^2-m^2n^2-12n^3+n^4}{(n+i)(n-i)(n-4)(n-5)(n+m)^2(n+m-2)}
\\
\times
{}_4F_3\left(\begin{array}{c}
-(n-m)/2+2,1+(n+m)/2-2,(m+1)/2,1+m/2\\
m+1,1+(m-i)/2,1+(m+i)/2
\end{array}\mid 1\right)
\\
-4
\frac{-96+64n+6m^2-m^2i^2+6i^2-14n^2-m^2n-ni^2+n^3)n(n-1)(n-m-2)(n-m-4)}
{(n+i)(n-i)(n-4)(n-6)(n+m)^2(n+m-2)(n+m-4)}
\\
\times
{}_4F_3\left(\begin{array}{c}
-(n-m)/2+3,1+(n+m)/2-3,(m+1)/2,1+m/2\\
m+1,1+(m-i)/2,1+(m+i)/2
\end{array}\mid 1\right)
\\
+
\frac{(n-m-4)(n-m-2)(n-2)(n-1)n(n-m-6)^2(n-6+i)(n-6-i)}
{(n+m-4)(n+m-2)(n+m)^2(n-6)(n-5)(n-4)(n-i)(n+i)}
\\
\times
{}_4F_3\left(\begin{array}{c}
-(n-m)/2+4,1+(n+m)/2-4,(m+1)/2,1+m/2\\
m+1,1+(m-i)/2,1+(m+i)/2
\end{array}\mid 1\right)
=0.
\end{multline}
\begin{rem}
Following a remark by Bailey \cite{BaileyPGMA2} it ought to be possible 
to reduce this 5-term relation to a 4-term relation.
\end{rem}
Switching $c_{n,m,i}$ back in with \eqref{eq.4F3}
and multiplication with the least common of all five denominators provides
\begin{multline}
(n+i)(n-i)(n-4)(n+m)(n-5)(n-6)(n-m)
c_{n,m,i}
\\
-
4(n-1)(4n-4n^2-m^2n-ni^2+n^3+m^2i^2)(n-5)(n-6)
c_{n-2,m,i}
\\
-2
(n-3)n(32-72n-2m^2-3m^2i^2-2i^2+48n^2+6m^2n-n^2i^2
\\
+6ni^2-m^2n^2-12n^3+n^4)(n-6)
c_{n-4,m,i}
\\
+4
(-96+64n+6m^2-m^2i^2+6i^2-14n^2-m^2n-ni^2+n^3)n(n-1)(n-5)
c_{n-6,m,i}
\\
+
(n-m-6)(n-2)(n-1)n(n-6+i)(n-6-i)(n+m-6)
c_{n-8,m,i}
=0.
\label{eq.c5T}
\end{multline}
\begin{rem}
The equation for $c_{n,m,i}$ in \eqref{eq.4F3} is symmetric with
respect to an exchange $(m-i)/2\leftrightarrow (m+i)/2$. 
This symmetry persists in \eqref{eq.c5T}. 
This is more obvious 
for example if we rewrite the second line
\begin{equation}
4n-4n^2-m^2n-ni^2+n^3+m^2i^2 = n(n-2)^2
-2n(m_+^2+m_-^2)
+(m_+^2-m_-^2)^2,
\end{equation}
or the third and fourth line
\begin{multline}
32-72n-2m^2-3m^2i^2-2i^2+48n^2+6m^2n-n^2i^2
+6ni^2-m^2n^2-12n^3+n^4
\\
=
(n-2)(n-4)(n^2-6n+4)
+(-4+12n-2n^2)(m_+^2+m_-^2)
-3(m_+^4+m_-^4)+6(m_+m_-)^2.
\end{multline}

\end{rem}

\subsection{The case $i\ge m$}
For $i$ in the range $m+2\ldots n$ small $s$-values
in the sum \eqref{eq.csum} are skipped:
\begin{multline}
c_{n,m,i}
=
(-)^{(n-m)/2}
\sum_{s=(i-m)/2}^{(n-m)/2} (-)^s
\binom{(n+m)/2+s}{(n-m)/2-s}\binom{m+2s}{s}
\\ \times
2^{1-m-2s}
{m+2s\choose (m+2s-i)/2}
\\
=
(-)^{(n-i)/2}
2^{1-i}
\binom{(n+i)/2}{i}\binom{i}{(i-m)/2}
{}_4F_3\left(
\begin{array}{c}
1+\frac{n+i}{2},
\frac{i+1}{2},
1+\frac{i}{2},
-\frac{n-i}{2}\\
1+\frac{m+i}{2},1+i,1+\frac{i-m}{2}
\end{array}
\mid 1
\right)
.
\label{eq.cnmiL}
\end{multline}

This $_4F_3$ has the again the structure covered by
\cite[(2.44)]{FieldsJAT1} and excess again 1/2, which means changing $n$ by 2 moves
two upper parameters contiguously up and down by 1. 
If one applies the series reversion formula
\cite{Mishevarxiv2020}
\begin{equation}
_{A+1}F_B[(a),-t; (b);z]
=
\frac{((a))_t(-z)^t}{((b))_t}
\,_{B+1}F_A[(1-b-t),-t;(1-a-t);\frac{(-1)^{A+B}}{z}]
\label{eq.rev}
\end{equation}
for positive integers $t$
to both polynomials, they turn out to have the common representation
\begin{multline}
c_{n,m,i} 
=
2^{1-n}
\binom{n}{(n-m)/2}
\binom{n}{(n-i)/2}
\\
\times
{}_4F_3\left(
\begin{array}{c}
-\frac{m+n}{2},
-\frac{i+n}{2},
-\frac{n-m}{2},
-\frac{n-i}{2}\\
-n,\frac{1-n}{2},-\frac{n}{2}
\end{array}
\mid 1
\right).
\label{eq.Falt}
\end{multline}
Because this is a symmetric function 
under the swap $i\leftrightarrow m$,
\begin{equation}
c_{n,m,i}=c_{n,i,m},
\label{eq.csym}
\end{equation}
it turns out that \eqref{eq.c5T}
also holds in the cases $i>m$, if the startup values \eqref{eq.m2m}
are replaced:
\begin{equation}
c_{n,m,i} 
= 
2^{1-i}
\binom{ (n+i)/2}{i}
\binom{i}{(i-m)/2}
\frac{\omega_{(n-i)/2}(i,m) }
{\prod_{d=1}^{(n-i)/2} (2d+i+m)(2d+i-m)},\quad i\ge m.
\label{eq.m2mhi}
\end{equation}
Note while using Table \ref{tab.omeg}
that the two arguments of $\omega$ in this formula are swapped in comparison
with the cases $i\le m$; the entries in the table are not 
symmetric under the involutions $m_-\leftrightarrow m_+$ or $i\leftrightarrow m$.

\section{Higher Dimensions}
Zernike radial polynomials in dimensions $D\ge 2$ are defined as
\begin{defn} (Radial Polynomials in Dimension $D$)
\begin{multline}
R_n^{(m)}(x)
=
(-1)^{(n-m)/2}
\sqrt{2n+D}
\binom{(D+n+m)/2-1}{(n-m)/2}
x^m
\\ \times
{}_2F_1\left(\begin{array}{c} -(n-m)/2, (D+n+m)/2 \\ m+D/2\end{array}\mid x^2\right).
\label{eq.RofxD}
\end{multline}
\end{defn}
This can also be written as Jacobi Polynomials
\begin{equation}
R_n^{(m)}(x)
=
(-1)^{(n-m)/2}
\sqrt{2n+D}
x^m
P_{(n-m)/2}^{m+D/2-1,0}(1-2x^2).
\end{equation}

The case $D=2$ is essentially the same as \eqref{eq.Rofx}, but we have
normalized them here as 
\begin{equation}
\int_0^1 R_n^{(m)}(x) R_{n'}^{(m)}(x) x^{D-1} dx =\delta_{n,n'}; \quad R_n^{(m)}(1)=\sqrt{2n+D},
\end{equation}
such that a $\sqrt{2n+D}$ factor appears on the right hand side
of the definition,  
and packed the upper parameter $m$ in parentheses
to indicate that difference. 
The factor $x^{D-1}$ in the integral is the radial part of the Jacobian determinant
of the hyperspherical coordinates.

If the difference between $n$ and $m$ is zero or two:
\begin{eqnarray}
R_n^{(n)}(x) & =& \sqrt{2n+D} x^n ; \label{eq.Rnn} \\
R_n^{(n-2)}(x) & =& -\sqrt{2n+D} 
\left(\frac{D}{2}+n-2\right)
x^{n-2} 
\left
[
1- \frac{D/2+n-1}{(n+D)/2-1} x^2
\right]
.
\end{eqnarray}

In lieu of Def.\ \ref{def.c} the coupling coefficients $c_{D,n,m,i}$ are defined as
\begin{equation}
R_n^{(m)}(x)\equiv \sqrt{2n+D} \mathop{{\sum}'}_i c_{D,n,m,i}T_i(x)
\label{eq.cdefD}
\end{equation}
for even $n-m\ge 0$. 
The square root has been factored out in that definition because 
\begin{itemize}
\item
 this keeps
\begin{equation}
c_{2,n,m,i} = c_{n,m,i}
\end{equation}
at $D=2$ aligned with $c_{n,m,i}$ of Section \ref{sec.M}, and 
\item
 the
$c_{D,n,m,i}$ are rational numbers that can be generated with software which has
support for exact representations of these---for reference purposes---see the appendix.
\end{itemize}

Inserting $x=1$ in \eqref{eq.cdefD} establishes the sum rule
\begin{equation}
\mathop{{\sum}'}_
{i=0 \atop n-i\,\mathrm{even}}^n
c_{D,n,m,i}=1.
\label{eq.srule}
\end{equation}

Expanding the  polynomial in \eqref{eq.RofxD} and replacing the monomials  with
\eqref{eq.xnofT} gives 
the generalized version of \eqref{eq.csum}:
\begin{multline}
c_{D,n,m,i}
=(-)^{(n-m)/2}
\sum_{s=\max(0,(i-m)/2)}^{(n-m)/2} 
(-)^s
\binom{(D+n+m)/2+s-1}{(n-m)/2} \binom{(n-m)/2}{s}
\\ \times
2^{1-m-2s}
{m+2s\choose (m+2s-i)/2}
\\
=
(-)^{(n-m)/2}2^{1-m}\binom{\frac{D+n+m}{2}-1}{\frac{n-m}{2}}\binom{m}{\frac{m-i}{2}}
{}_4F_3\left(\begin{array}{c}
\frac{D+n+m}{2}, -\frac{n-m}{2},\frac{m+1}{2}, \frac{m}{2}+1 \\
\frac{D}{2}+m, 1+\frac{m-i}{2}, 1+\frac{m+i}{2}
\end{array}\mid 1\right)
\label{eq.4F3D}
\end{multline}
if $0\le i \le n$ and $0\le m \le n$ and $2 \mid (n-m)$ and $2\mid (n-i)$, otherwise $c_{D,n,m,i}$ is zero.
The excess of this $_4F_3$ is again 1/2.

Contiguous relations of the Gaussian Hypergeometric Functions or
the equivalent recurrences for Jacobi Polynomials induce
the 3-term recurrence  \cite{LewanowiczZMAM17,MatharArxiv0705a} 
\begin{multline}
-(1+\frac{n-m}{2})(1-n-\frac{D}{2})\frac{n+m+D}{2}
\frac{R_{n+2}^{(m)}(x)}{\sqrt{2(n+2)+D}}
\\
=
\frac{n-m}{2}(1+n+\frac{D}{2})(1-\frac{n+m+D}{2})
\frac{R_{n-2}^{(m)}(x)}{\sqrt{2(n-2)+D}}
\\
+(n+\frac{D}{2})\left[(1+n+\frac{D}{2})(1-n-\frac{D}{2})(1-x^2) 
+\frac12(n-m)(D+n+m)
+m+\frac{D}{2}-1\right]
\frac{R_n^{(m)}(x)}{\sqrt{2n+D}}
.
\end{multline} 
For each of the three $R_n^{(m)}$ the expansion \eqref{eq.cdefD} is inserted. $1-x^2$
is replaced by $[1-T_2(x)]/2$ \cite[Tab. 22.3]{AS} and  \cite[22.7.24]{AS}
\begin{equation}
T_i(x)T_j(x)=\frac12[T_{|i-j|}(x)+T_{i+j}(x)]
\end{equation}
is used on the right hand side
to replace products of Chebyshev Polynomials by  single terms.
Then both sides are projected onto a single $T_i$ by multiplying them
with $T_i(x)/\sqrt{1-x^2}$ and integrating using the
orthogonality \eqref{eq.Tortho}.

\begin{multline}
-(1+\frac{n-m}{2})(1-n-\frac{D}{2})\frac{n+m+D}{2}
c_{D,n+2,m,i}
\\
=
\frac{n-m}{2}(1+n+\frac{D}{2})(1-\frac{n+m+D}{2})
c_{D,n-2,m,i}
\\
-\frac14(n+\frac{D}{2})
(1+n+\frac{D}{2})(1-n-\frac{D}{2})
\\ \times
\left[
\frac12 c_{D,n,m,0} \delta_{2,i}
+
c_{D,n,m,1} \delta_{1,i}
+
c_{D,n,m,i+2} \epsilon_i
+c_{D,n,m,i-2} \frac{[i\ge 2]}{\epsilon_{i-2}}
\right]
\\
+(n+\frac{D}{2})
\left[
\frac12
(1+n+\frac{D}{2})(1-n-\frac{D}{2})
+
\frac12(n-m)(D+n+m)
+m+\frac{D}{2}-1\right]
c_{D,n,m,i}
\label{eq.recD}
\end{multline} 
for $n\ge 0$.
The notation $[.]$ in the numerator in the fourth line is the Iverson bracket: 
the square bracket equals 0 if the condition in the bracket is false,
else 1.
The recurrence starts at
\begin{equation}
c_{D,0,0,0} = 2;\quad c_{D,1,1,1} =1
\label{eq.D3init}
\end{equation}
and assumes that $c$ is zero whenever the requirements on even parities of $n-m$ and $n-i$
or $0\le m,i \le n$ are not fulfilled.

Examples are in Table \ref{tab.cD}.
Note that the symmetry \eqref{eq.csym} does not hold for $D\neq 2$.
Equivalent to \eqref{eq.Rnn} a restatement of \eqref{eq.xnofT} is
\begin{equation}
c_{D,n,n,i}=
2^{1-n}
\binom{n}{(n-i)/2},
\label{eq.mm}
\end{equation}
so the $c$-values 
in Table \ref{tab.cD} 
where $n=m$ 
do not depend on $D$. More generally, $c_{D,n,m,i}$ is a polynomial of order $\kappa$
in $D$---consequence of the fact that the $D$-dependent
binomial factor in the first line of \eqref{eq.4F3D} has $\kappa$ as the lower parameter.

\begin{rem}\label{rem.stor}
The growth of this list of constants for some fixed $D$ is as follows: for fixed $n$ and $m$
there are $\lfloor 1+n/2\rfloor$ entries---one for each $i$ that matches
the parity requirement. For fixed $n$ and variable $m$ and $i$ there are 
$(\lfloor 1+n/2\rfloor)^2 = [2n^2+6n+5+(2n+3)(-1)^n]/8$ entries. The number of entries for $n$ up
to some upper limit $\hat n$ is the sum over this over $n=0,\ldots, \hat n$,
$(\hat n+2)[2\hat n^2+8\hat n+9+3(-1)^{\hat n}]/24$.
This growth roughly 
$\propto
\hat n^3/12$ is the main obstacle of implementing this in memory-aware
programs to compute $R_n^{(m)}(x)$.
[For $D=2$ the space needed to store the coefficients is halved if the symmetry \eqref{eq.csym} is exploited.]
\end{rem}

\begin{table}
\begin{tabular}{rr|cccccc}
$n$ & $m$ & $c_{D,n,m,i}$ \\
\hline
0& 0& 2 \\
1& 1& 1 \\
2& 0& $1-\frac{D}{2}$ & $\frac12+\frac{D}{4}$ \\
2& 2& 1& $\frac12$ \\
3& 1& $\frac12-\frac{D}{8}$ & $\frac12+\frac{D}{8}$ \\
3& 3& 3/4& 1/4 \\
4& 0& $\frac14-\frac{D}{16}+\frac{3D^2}{32}$ & $-\frac{(D-2)(D+4)}{16}$ & $\frac{(D+4)(D+6)}{64}$ \\
4& 2& $\frac14-\frac{D}{8}$ & $\frac12$ & $\frac38+\frac{D}{16}$ \\
4& 4& $\frac34$& $\frac12$ & $\frac18$ \\
5& 1& $\frac14-\frac{D}{32}+\frac{D^2}{64}$& $-\frac{(D+6)(3D-8)}{128}$ 
   & $\frac{(D+6)(D+8)}{128}$ \\
5& 3& $\frac14-\frac{D}{16}$ & $\frac12+\frac{D}{32}$ & $\frac14+\frac{D}{32}$ \\
5& 5& $\frac58$ & $\frac{5}{16}$ & $\frac{1}{16}$ \\
6& 0& $-\frac{(D-2)(5D^2+22D+48)}{384}$ & $\frac{(D+6)(5D^2-6D+16)}{512}$ 
   & $-\frac{(D-2)(D+6)(D+8)}{256}$ 
   & $\frac{(D+6)(D+8)(D+10)}{1536}$ \\
6& 2& $\frac14+\frac{D}{32}+\frac{D^2}{64}$ & $\frac{3}{16}-\frac{9D}{128}-\frac{D^2}{256}$ 
   & $-\frac{(D-6)(D+8)}{128}$ 
   & $\frac{(D+8)(D+10)}{256}$ \\
6& 4& $\frac18-\frac{D}{16}$ & $\frac{11}{32}-\frac{D}{64}$ 
   & $\frac{7}{16}+\frac{D}{32}$ 
   & $\frac{5}{32}+\frac{D}{64}$ \\
6& 6& $\frac58$ & $\frac{15}{32}$ 
   & $\frac{3}{16}$ 
   & $\frac{1}{32}$ \\
\end{tabular}
\caption{Reference values of $c_{D,n,m,i}$ complete up to $n\le 6$. The $c$-values 
are sorted left-to-right by increasing $i=0,2,4,\ldots$ if $n$ is even, by
increasing $i=1,3,5,7,\ldots$ if $n$ is odd. The case $D=2$
turns this into Table \ref{tab.c}.
}
\label{tab.cD}
\end{table}

The recurrence \eqref{eq.recD} differs in spirit from the recurrence in Section \ref{sec.M} (besides
covering the wider range $D\ge 2$):
here three different $n$-values are coupled to three different $i$-values
at constant $m$, whereas in Section \ref{sec.M} five different $n$-values are
coupled at constant $m$ and $i$.

Rewriting \eqref{eq.4F3D} with the term reversion \eqref{eq.rev} gives
\begin{multline}
c_{D,n,m,i}
=
2^{1-n}
\binom{n}{(n-i)/2} \binom{n+D/2-1}{(n-m)/2)}
\\ \times
{}_4F_3\left(\begin{array}{c}
1-\frac{D+n+m}{2}, -\frac{n-m}{2},-\frac{n-i}{2}, -\frac{n+i}{2} \\
-\frac{n-1}{2}, 1-n-\frac{D}{2},-\frac{n}{2}
\end{array}\mid 1\right)
,
\end{multline}
the generalized \eqref{eq.Falt}.
In \eqref{eq.4F3D} modifying $n$
in steps of 2 modifies two of the upper parameters in the $_4F_3$
contiguously in steps  of 1, keeping the other five parameters constant.
\emph{Here}, modifying $m$ or $i$
in steps of 2 modifies two of the upper parameters in the $_4F_3$
contiguously in steps  of 1, keeping the other five parameters constant.
In both cases, the recurrence \eqref{eq.5rec} applies:
5-term recurrences exist that couple (i) 5
different $n$ at constant $\{D, m, i\}$,
or (ii) couple 5 different $m$
at constant $\{D, n, i\}$, 
or (iii) couple 5 different $i$
at constant $\{D, n, m\}$. 
There may be an advantage of such an algorithm compared with \eqref{eq.recD} 
in a threaded
program because the input parameters are better disentangled.

We present the explicit recurrence only for algorithm (iii):
\begin{multline}
(n+i+D)(i+5)(n+i+2)(n+i+4)(n+i+6)
\frac{n-i}{2}
(\frac{n-i}{2}-1)
(\frac{n-i}{2}-2)
(\frac{n-i}{2}-3)
c_{D,n,m,i}
\\
+2(n+i+4)(n+i+6)
(\frac{n-i}{2}-1)
(\frac{n-i}{2}-2)
(\frac{n-i}{2}-3)
(\frac{n+i}{2}+1)
\\ \times
(-24-18i+8D-3i^2+n^2+nD+i^2D+6iD-2Dmi-10Dm+4im+20m-2m^2i-10m^2)
c_{D,n,m,i+2}
\\
-2(i+4)(n-i-2)(n+i+6)
(\frac{n-i}{2}-2)
(\frac{n-i}{2}-3)
(\frac{n+i}{2}+1)
(\frac{n+i}{2}+2)
\\ \times
(-12-8i-D-i^2+n^2+nD+8m-4Dm-4m^2)
c_{D,n,m,i+4}
\\
-2(n-i-2)(n-i-4)
(\frac{n-i}{2}-3)
(\frac{n+i}{2}+1)
(\frac{n+i}{2}+2)
(\frac{n+i}{2}+3)
\\ \times
(-72-30i+24D-3i^2+n^2+nD+i^2D+10iD-4im+2Dmi+2m^2i-12m+6Dm+6m^2)
c_{D,n,m,i+6}
\\
+(n-i-2)(n-i-4)(n-i-6)
(\frac{n+i}{2}+1)
(\frac{n+i}{2}+2)
(\frac{n+i}{2}+3)
(\frac{n+i}{2}+4)
\\ \times
(n-i-8+D)(i+3)
c_{D,n,m,i+8}
=0.
\end{multline}

\section{Summary}
The coupling coefficients $c_{n,m,i}$  and $c_{D,n,m,i}$
are nonzero
only if $0\le i \le n$, $0\le m\le n$ and if $m-i$ and $n-i$ are even.
\subsection{Circle Polynomials}
A calculation for $D=2$ may start from 
\begin{itemize}
\item
the product expressions
in \eqref{eq.m2m} if $0\le (n-m)/2\le 3$  and $i\le m$ 
and the aid of Table \ref{tab.omeg}
\item
the cross-reference \eqref{eq.csym} if $m\le i$,
\end{itemize}
and calculate the values
for larger $n-m$ with the 5-term recurrence \eqref{eq.c5T}.

If $i=0$, \eqref{eq.cnm0} may be used. If $i=1$, the 3-term recurrence
\eqref{eq.rec1} is also convenient.

\subsection{General $D$}
A generic calculation for $D\ge 2$ may start from  \eqref{eq.D3init}
and calculate the others with the 5-term recurrence \eqref{eq.recD}.

\appendix
\section{Software}
The arXiv ancillary directory 
\texttt{anc/cnmi.asc}
contains Python, Maple, Pari/gp, Mathematica and C++ programs
to produce the table in the file \texttt{cnmi.asc}. Each line in the table 
\texttt{anc/cnmi.asc} shows the four integer
values $D$, $n$, $m$ and $i$, followed by a floating point and an exact representation
of $c_{D,n,m,i}$.
These programming languages have multi-precision capabilities (integer
libraries) to produce exact
results even though the alternating fractions in the sum over the $s$ may have large binomial products
and may suffer from cancellations of digits in fixed-precision representations of numbers.
\begin{rem}
All these implementations perform the summations in \eqref{eq.4F3D}
over a number of terms that grows $\propto n$;
they are obviously not optimized for speed.
\end{rem}

\begin{rem}
The Pari/gp script needs 
less than a second to produce that table for $D=2$ and $D=3$ for $n\le 40$
on a i7-2600 3.40 GHz CPU,
and 3 minutes to produce a complete table of $c_{2,n,m,i}$ for $n\le 200$.
\end{rem}

\bibliographystyle{amsplain}
\bibliography{all}

\end{document}